\documentclass[11pt]{article}
\usepackage{amsfonts, amssymb, amsmath, amsthm, xcolor, xy}
\xyoption{arc}
\begin{document}
\newcommand{\m}{\textnormal{mult}}
\newcommand{\ra}{\rightarrow}
\newcommand{\la}{\leftarrow}
\renewcommand{\baselinestretch}{1.1}

\theoremstyle{plain}
\newtheorem{thm}{Theorem}[section]
\newtheorem{cor}[thm]{Corollary}
\newtheorem{con}[thm]{Conjecture}
\newtheorem{cla}[thm]{Claim}
\newtheorem{lm}[thm]{Lemma}
\newtheorem{prop}[thm]{Proposition}
\newtheorem{example}[thm]{Example}

\theoremstyle{definition}
\newtheorem{dfn}[thm]{Definition}
\newtheorem{alg}[thm]{Algorithm}
\newtheorem{rem}[thm]{Remark}

\renewcommand{\baselinestretch}{1.1}

\title{\bf Maximum Multiplicity of a Root of the Matching Polynomial of a Tree and Minimum Path Cover}
\author{
Cheng Yeaw Ku
\thanks{ Department of Mathematics, National University of Singapore, Singapore 117543. E-mail: matkcy@nus.edu.sg} \and K.B. Wong \thanks{
Institute of Mathematical Sciences, University of Malaya, 50603 Kuala Lumpur, Malaysia. E-mail:
kbwong@um.edu.my.} } \maketitle

\begin{abstract}\noindent
We give a necessary and sufficient condition for the maximum multiplicity of a root of the matching polynomial of a tree to be equal to the minimum number of vertex disjoint paths needed to cover it.
\end{abstract}

\section{Introduction}

All the graphs in this paper are simple. The vertex set and the edge set of a graph $G$ are denoted by $V(G)$ and $E(G)$ respectively. A {\em matching} of a graph $G$ is a set of pairwise disjoint edges of $G$. Recall that for a graph $G$ on $n$ vertices, the {\it matching
polynomial} $\mu(G,x)$ of $G$ is given by
\[ \mu(G,x)=\sum_{k \ge 0} (-1)^{k}p(G,k)x^{n-2k}, \]
where $p(G,k)$ is the number of matchings with $k$ edges in $G$. Let $\m(\theta, G)$ denote the multiplicity
of $\theta$ as a root of $\mu(G,x)$.

The following results are well known. The proofs can be found in \cite[Theorem 4.5 on p. 102]{G0}.

\begin{thm}\label{max-path}
The maximum multiplicity of a root of the matching polynomial $\mu(G, x$) is at most the minimum number of vertex disjoint paths needed to cover the vertex set of $G$.
\end{thm}

Consequently,

\begin{thm}\label{path}
If $G$ has a Hamiltonian path, then all roots of its matching polynomial are simple.
\end{thm}

The above is the source of motivation for our work. It is natural to ask when does equality holds in Theorem \ref{max-path}. In this note, we give a necessary and sufficient condition for the maximum multiplicity of a root of the matching polynomial of a tree to be equal to the minimum number of vertex disjoint paths needed to cover it. Before stating the main result, we require some terminology and basic properties of matching polynomials.

It is well known that the roots of the matching polynomial are real. If $u\in V(G)$, then $G \setminus u$ is the graph obtained from $G$ by deleting the vertex $u$ and the edges of $G$ incident to $u$. It is known that the roots of $G\setminus u$ {\em interlace} those of $G$, that is, the multiplicity of a root changes by at most one upon deleting a vertex from $G$. We refer the reader to \cite{G0} for an introduction to matching polynomials.

\begin{lm}\label{interlacing}
Suppose $\theta$ is a root of $\mu(G,x)$ and $u$ is a vertex of $G$. Then
\[ \m(\theta, G)-1 \le \m(\theta, G \setminus u) \le \m(\theta, G)+1. \]
\end{lm}

\noindent As a consequence of Lemma \ref{interlacing}, we can classify the vertices in a graph by assigning a `sign' to each vertex (see \cite{G}).

\begin{dfn}
Let $\theta$ be a root of $\mu(G,x)$. For any vertex $u\in V(G)$,
\begin{itemize}
\item $u$ is $\theta$-{\em essential} if $\m(\theta, G
\setminus u)=\m(\theta, G)-1$,

\item $u$ is $\theta$-{\em neutral} if $\m(\theta, G
\setminus u)=\m(\theta, G)$,

\item $u$ is $\theta$-{\em positive} if $\m(\theta, G
\setminus u)=\m(\theta, G)+1$.
\end{itemize}
\end{dfn}

Clearly, if $\m(\theta, G)=0$ then there are no $\theta$-essential
vertices since the multiplicity of a root cannot be negative.
Nevertheless, it still makes sense to talk about $\theta$-neutral
and $\theta$-positive vertices when $\m(\theta, G)=0$. The converse
is also true, i.e. any graph $G$ with $\m(\theta, G)>0$ must have at
least one $\theta$-essential vertex. This was proved in \cite[Lemma
3.1]{G}.

\noindent A further classification of vertices plays an important role in establishing some structural properties of a graph:
\begin{dfn}Let $\theta$ be a root of $\mu(G,x)$. For any vertex $u\in V(G)$, $u$ is $\theta$-{\em special} if it is
not $\theta$-essential but has a neighbor that is $\theta$-essential.
\end{dfn}

\noindent If $G$ is connected and not all of its vertices are $\theta$-essential, then $G$ must contain a $\theta$-special vertex. It turns out that a $\theta$-special vertex must be $\theta$-positive (see \cite[Corollary 4.3]{G}).

We now introduce the following definition which is crucial in describing our main result.

\begin{dfn}
Let $G$ be a graph and $\mathcal{Q}=\{Q_1,\ldots, Q_m\}$ be a set of vertex disjoint paths that cover $G$. Then $\mathcal{Q}$ is said to be $(\theta, G)$-{\em extremal} if it satisfies the following:
\begin {itemize}
\item [(a)] $\theta$ is a root of $\mu (Q_i,x)$ for all $i=1, \ldots, m$;
\item [(b)] for every edge $e=\{u, v\} \in E(G)$ with $u \in Q_r$ and $v \in Q_s$, $r \neq s$, either
 $u$ is $\theta$-special in $Q_r$ or $v$ is $\theta$-special in $Q_s$.
\end {itemize}
\end{dfn}

Our main result is the following:

\begin{thm}\label{main}
Let $T$ be a graph and $\mathcal{Q}=\{Q_1, \ldots, Q_m\}$ be a set of vertex disjoint paths covering $T$. Then $m$ is the maximum multiplicity of a root of the matching polynomial $\mu(T,x)$, say $\m(\theta, T)=m$ for some root $\theta$, if and only if $\mathcal{Q}$ is $(\theta, T)$-extremal.
\end{thm}

The following example shows that Theorem \ref{main} cannot be extended to general graphs.

\begin{example}
Consider the following graph $G$:
\end{example}

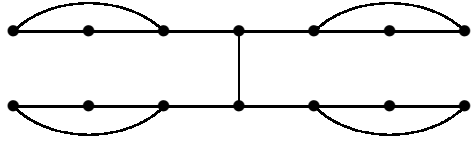
\[\xy
(0,0)*{}="A"; (10,0)*{}="B"; (20,0)*{}="C";
(30,0)*{}="D";(40,0)*{}="E"; (50,0)*{}="F"; (60,0)*{}="G";
"A";"B"**\dir{-};"B";"C"**\dir{-};"C";"D"**\dir{-};"D";"E"**\dir{-};"E";"F"**\dir{-};"F";"G"**\dir{-};
(0,0)*{\bullet};(10,0)*{\bullet};(20,0)*{\bullet};(30,0)*{\bullet};(40,0)*{\bullet};(50,0)*{\bullet};(60,0)*{\bullet};
(0,-10)*{}="A2"; (10,-10)*{}="B2"; (20,-10)*{}="C2";
(30,-10)*{}="D2";(40,-10)*{}="E2"; (50,-10)*{}="F2"; (60,-10)*{}="G2";
"A2";"B2"**\dir{-};"B2";"C2"**\dir{-};"C2";"D2"**\dir{-};"D2";"E2"**\dir{-};"E2";"F2"**\dir{-};"F2";"G2"**\dir{-};
(0,-10)*{\bullet};(10,-10)*{\bullet};(20,-10)*{\bullet};(30,-10)*{\bullet};(40,-10)*{\bullet};(50,-10)*{\bullet};(60,-10)*{\bullet};
"D";"D2"**\dir{-}; "A";"C"**\crv{(5,5) & (15,5)};"E";"G"**\crv{(45,5) & (55,5)};"D";"D2"**\dir{-}; "A2";"C2"**\crv{(5,-15) & (15,-15)};"E2";"G2"**\crv{(45,-15) & (55,-15)};
\endxy\]

Let $P_{7}$ denote the path on $7$ vertices. Note that $\m(\sqrt{3}, G) = 2$ and $\mu(P_{7}, x) = x^{7}-6x^{5}+10x^{3}-4x$. By Theorem \ref{max-path}, the maximum multiplicity of a root of $\mu(G,x)$ is $2$. Also, $G$ can be covered by two paths on $7$ vertices. However, $\sqrt{3}$ is not a root of $\mu(P_{7}, x)$.

\section{Basic Properties}

In this section, we collect some useful results proved in \cite{G0} and \cite{G}. Recall that if $u\in V(G)$, then $G \setminus u$ is the graph obtained from $G$ by deleting vertex $u$ and the edges of $G$ incident to $u$. We also denote the graph $(G \setminus u) \setminus v$ by $G \setminus uv$. Note that the resulting graph does not depend on the order of which the vertices are deleted.

If $e\in E(G)$, the graph $G-e$ is the graph obtained from $G$ by deleting the edge $e$.

The matching polynomial satisfies the following basic identities, see \cite[Theorem 1.1 on p. 2]{G0}.
\begin{prop}\label{identity}
Let $G$ and $H$ be graphs, with matching polynomials $\mu(G,x)$ and $\mu(H,x)$, respectively. Then
\begin {itemize}
\item [\textnormal{(a)}] $\mu(G \cup H,x) = \mu(G,x)\mu(H,x)$,
\item [\textnormal{(b)}] $\mu(G,x)=\mu(G - e, x)-\mu(G
\setminus uv,x)$ where $e=\{u,v\}$ is an edge of $G$,
\item [\textnormal{(c)}] $\mu(G,x) = x\mu(G \setminus u,x)-\sum_{v
\sim u} \mu(G \setminus uv,x)$ for any vertex $u$ of $G$.
\end {itemize}
\end{prop}

Suppose $P$ is a path in $G$. Let $G \setminus P$ denote the graph
obtained from $G$ by deleting the vertices of $P$ and all the edges
incident to these vertices. It is known that the multiplicity of a root decreases by at most one upon deleting a path, see \cite[Corollary 2.5]{G}.

\begin{lm}\label{interlacing-path}
For any root $\theta$ of $\mu(G,x)$ and a path $P$ in $G$,
\[ \m(\theta, G \setminus P) \ge \m(\theta, G)-1. \]
\end{lm}

\noindent If equality holds, we say that the path $P$ is $\theta$-{\em essential} in $G$. Godsil \cite{G} proved that if a vertex $v$ is not $\theta$-essential in $G$, then no path with $v$ as an end point is $\theta$-essential. In other words,

\begin{lm}\label{essential-path}
If $P$ is a $\theta$-essential path in $G$, then its endpoints are $\theta$-essential in $G$.
\end{lm}

The next result of Godsil \cite[Corollary 4.3]{G} implies that a $\theta$-special vertex must be $\theta$-positive.

\begin{lm}\label{neutral-essential}
A $\theta$-neutral vertex cannot be joined to any $\theta$-essential vertex.
\end{lm}

\section{Gallai-Edmonds Decomposition}

It turns out that $\theta$-special
vertices play an important role in the Gallai-Edmonds decomposition
of a graph. We now define such a decomposition. For any root
$\theta$ of $\mu(G,x)$, partition the vertex set $V(G)$ as follows:
\begin{eqnarray*}
D_{\theta}(G) & = & \{ u : u ~~\rm{is }~~\theta\textnormal{-essential in}~ G \} \\
A_{\theta}(G) & = & \{ u : u~~\rm{is }~~\theta\textnormal{-special in}~ G \}\\
C_{\theta}(G) & = & V(G)-D_{\theta}(G)-A_{\theta}(G).
\end{eqnarray*}
We call these sets of vertices the $\theta$-{\em partition classes}
of $G$. The Gallai-Edmonds Structure Theorem is usually stated in
terms of the structure of maximum matchings of a graph with respect
to its $\theta$-partition classes when $\theta=0$. Its proof
essentially follows from the following assertions (for more
information, see \cite[Section 3.2]{LP}):

\begin{thm}[Gallai-Edmonds Structure Theorem]\label{GE}\hfill

\noindent Let $G$ be any graph and let $D_{0}(G)$, $A_{0}(G)$ and
$C_{0}(G)$ be the $0$-partition classes of $G$.
\begin{itemize}
\item[\rm{(i)}] \textnormal{({\em The Stability Lemma})} Let $u
\in A_{0}(G)$ be a $0$-special vertex in $G$. Then

$\bullet$ $v \in D_{0}(G)$ if and only if $v \in D_{0}(G \setminus
u)$;

$\bullet$ $v \in A_{0}(G)$ if and only if $v \in A_{0}(G \setminus
u)$;

$\bullet$ $v \in C_{0}(G)$ if and only if $v \in C_{0}(G \setminus
u)$.

\item[\rm{(ii)}] \textnormal{({\em Gallai's Lemma})} If every
vertex of $G$ is $0$-essential then $\m(0, G)=1$.
\end{itemize}
\end{thm}
For any root $\theta$ of $\mu(G,x)$, it was shown by Neumaier
\cite[Corollary 3.3]{N} that the analogue of Gallai's Lemma holds when $G$ is a
tree. A different proof was given by Godsil (see \cite[Corollary
3.6]{G}).

\begin{thm}[\cite{G}, \cite{N}]\label{Neu}
Let $T$ be a tree and let $\theta$ be a root of $\mu(T,x)$. If every
vertex of $T$ is $\theta$-essential then $\m(\theta, G)=1$.
\end{thm}

On the other hand, it was proved in \cite[Theorem 5.3]{G} that if
$\theta$ is any root of $\mu(T,x)$ where $T$ is tree and $u \not \in
D_{\theta}(T)$, then $v \in D_{\theta}(T)$ if and only if $v \in
D_{\theta}(T \setminus u)$. It turns out that this assertion is
incorrect (see Example \ref{counter-example} below). However, using the idea of the proof of Theorem 5.3 in \cite{G}, we shall prove the Stability Lemma for trees with any given root of its
matching poynomial. Note that the Stability Lemma is a weaker statement than Theorem 5.3 in \cite{G}. Together with Theorem \ref{Neu}, this yields
the Gallai-Edmonds Structure Theorem for trees with general root
$\theta$. Recently, Chen and Ku \cite{KC} had proved the Gallai-Edmonds Structure Theorem for general graph with any root $\theta$. However, our proof of the special case for trees, which uses an eigenvector argument, is different from the the one given in \cite{KC}. We believe that different proofs can be illuminating. For the sake of completeness, we include the proof in the next section.

\begin{thm}[The Stability Lemma for Trees]\label{stability-tree}
Let $T$ be a tree and let $\theta$ be a root of $\mu(T,x)$. Let $u
\in A_{\theta}(T)$ be a $\theta$-special vertex in $T$. Then

$\bullet$ $v \in D_{\theta}(T)$ if and only if $v \in D_{\theta}(T
\setminus u)$;

$\bullet$ $v \in A_{\theta}(T)$ if and only if $v \in A_{\theta}(T
\setminus u)$;

$\bullet$ $v \in C_{\theta}(T)$ if and only if $v \in C_{\theta}(T
\setminus u)$.
\end{thm}

It is well known that the matching polynomial of a graph $G$ is equal to the characteristic polynomial of $G$ if and only if $G$ is a forest. To prove Theorem \ref{stability-tree}, the following characterization of $\theta$-essential vertices in a tree via eigenvectors is very useful. Recall that a vector $f \in \mathbb{R}^{|V(G)|}$ is an eigenvector of a graph $G$ with eigenvalue $\theta$ if and only if for every vertex $u \in V(G)$,
\begin{eqnarray}
\theta f(u) = \sum_{v \sim u} f(v). \label{eigenvalue}
\end{eqnarray}

\begin{prop}[{\cite[Theorem 3.4]{N}}]\label{Neumaier}
Let $T$ be a tree and let $\theta$ be a root of its matching
polynomial. Then a vertex $u$ is $\theta$-essential if and only if
there is an eigenvector $f$ of $T$ such that $f(u) \not = 0$.
\end{prop}

\noindent In fact Proposition \ref{Neumaier} can be deduced from Lemma 5.1 of \cite{G}. The following corollaries are immediate consequences of
Proposition \ref{Neumaier}.

\begin{cor}\label{C1}
Let $T$ be a tree and let $\theta$ be a root of its matching
polynomial. Let $\alpha \not = 0$ be a nonzero real number. If $f$
is an eigenvector of $T$ such that $f(u) \not = 0$, then there
exists another eigenvector $g$ such that $g(u)=\alpha \not = 0$.
Moreover, both $g$ and $f$ have the same support, i.e. $\{i: f(i)
\not = 0\}=\{i: g(i) \not = 0\}$.
\end{cor}

\begin{proof} Let $g \in \mathbb{R}^{|V(T)|}$ be the vector defined as follows: for any vertex $x$ of
$T$,
\[ g(x) = \frac{\alpha}{f(u)} \cdot f(x). \]
Then, for any vertex $x$ of $T$,
\begin{eqnarray}
\sum_{y \sim x} g(y) & = & \sum_{y \sim x} \frac{\alpha}{f(u)} \cdot f(y) \nonumber \\
& = & \frac{\alpha}{f(u)} \sum_{y \sim x} f(y) \nonumber \\
& = & \frac{\alpha}{f(u)} \theta f(x) ~~\textnormal{(as $f$ is an eigenvector)}\nonumber \\
& = & \theta g(x). \nonumber
\end{eqnarray}
So $g$ is the required vector. Clearly, $g(x)=0$ if and only if
$f(x)=0$, so $g$ and $f$ have the same support.
\end{proof}

\begin{cor}[{\cite[Theorem 3.4]{N}}]\label{C2}
Let $T$ be a tree. If $u$ is $\theta$-special then it is joined to
at least two $\theta$-essential vertices.
\end{cor}

\begin{proof} By definition, $u$ has a $\theta$-essential neighbor, say
$w$. By Proposition \ref{Neumaier}, there exists an eigenvector $f$
of $T$ corresponding to $\theta$ such that $f(v) \not = 0$. By
Proposition \ref{Neumaier} again, $f(u)=0$ and so $\sum_{v \sim
u}f(v)=0$. This implies that $f(v) \not = 0$ on at least two
neighbors of $u$. By Proposition \ref{Neumaier}, both of them are
$\theta$-essential.
\end{proof}

The following assertion follows from Theorem \ref{Neu} and Proposition \ref{Neumaier}.

\begin{cor}[{\cite[Corollary 3.3]{N}}]\label{C3}
Let $T$ be a tree and let $\theta$ be a root of $\mu(T,x)$. Suppose every vertex of $T$ is $\theta$-essential. Then every $\theta$-eigenvector of $T$ has no zero entries.
\end{cor}

We also require the following partial analogue of the Stability
Lemma for general root obtained by Godsil in \cite{G}.

\begin{prop}[{\cite[Theorem 4.2]{G}}]\label{delete-positive}
Let $\theta$ be a root of $\mu(G,x)$ with non-zero multiplicity $k$
and let $u$ be a $\theta$-positive vertex in $G$. Then
\begin{itemize}
\item[\textnormal{(a)}] if $v$ is $\theta$-essential in $G$ then it
is $\theta$-essential in $G \setminus u$;
\item[\textnormal{(b)}] if $v$ is $\theta$-positive in $G$ then it
is $\theta$-essential or $\theta$-positive in $G \setminus u$;
\item[\textnormal{(c)}] if $u$ is $\theta$-neutral in $G$ then it is
$\theta$-essential or $\theta$-neutral in $G \setminus u$.
\end{itemize}
\end{prop}

\section{Proof of the Stability Lemma for Tree}

This section is devoted to the proof of Theorem \ref{stability-tree}, which will follow from the following theorem.

\begin{thm}\label{all}
Let $T$ be a tree and let $\theta$ be a root of $\mu(T,x)$. Then there exists a $\theta$-eigenvector $f$ of $T$ such that $f(x) \not = 0$ for every $\theta$-essential vertex $x$ in $T$. Moreover, if $v$ is $\theta$-essential in $T \setminus u$ where $u$ is $\theta$-special in $T$, then $v$ is $\theta$-essential in $T$.
\end{thm}

\begin{proof}
If every vertex of $T$ is $\theta$-essential, then the result follows from Corollary \ref{C3}. Therefore, we may assume that $T$ has a $\theta$-special vertex, say $u$. We proceed by induction on the number of vertices.

Suppose $b_{1}, \ldots, b_{s}$ are all the neighbors of $u$ in $T$. Then each $b_{i}$ belongs to different components of $T \setminus u$, say $b_{i} \in V(C_{i})$ where $C_{1}, \ldots C_{s}$ are components of $T \setminus u$.

First, we partition the set $\{1, \ldots, s\}$ as follows:

\begin{eqnarray}
A & = & \{i: b_{i}~~\textrm{is } \textrm{$\theta$-essential in }C_{i}\}, \nonumber \\
B & = & \{ i: i \not \in A, \theta \textrm{ is a root of } \mu(C_{i}, x)\}, \nonumber \\
C & = & \{1, \ldots, s\} \setminus (A \cup B). \nonumber
\end{eqnarray}

By the inductive hypothesis, for each $i \in A$, there exists a $\theta$-eigenvector $f_{i}$ of $C_{i}$ such that $f_{i}(x) \not = 0$ for every $\theta$-essential vertex $x$ in $C_{i}$. In particular, $f_{i}(b_{i}) \not = 0$ for all $i \in A$. By Proposition \ref{delete-positive}, any $\theta$-essential vertex in $T$ is also $\theta$-essential in $T \setminus u$, so for each $i \in A$,
\begin{eqnarray}
f_{i}(x) \not = 0~~~~\textrm{if } x ~\textrm{is $\theta$-essential in $T$ and } x \in V(C_{i}). \label{e2}
\end{eqnarray}
For every $i \in A$, choose $\alpha_{i} \in \mathbb{R}$ such that
\[ \alpha_{i} \not = 0 ~~~~\textrm{and}~~~~\sum_{i \in A} \alpha_{i} = 0. \]
Such a choice is always possible since $|A| \ge 2$ by Corollary \ref{C2}. Now, applying Corollary \ref{C1} to each $f_{i}$ with $i \in A$, we obtain the eigenvectors $g_{i}$ of $C_{i}$ such that $g_{i}(b_{i}) = \alpha_{i} \not = 0$ with both $g_{i}$ and $f_{i}$ having the same support. In particular, it follows from (\ref{e2}) that for each $i \in A$,
\begin{eqnarray}
g_{i}(x) \not = 0~~~~\textrm{if } x ~\textrm{is $\theta$-essential in $T$ and } x \in V(C_{i}). \label{e3}
\end{eqnarray}

Also, for each $i \in B$, by the inductive hypothesis, we can choose an eigenvector $g_{i}$ so that $g_{i}(x) \not = 0 $ for every $\theta$-essential vertex $x$ in $C_{i}$. By Proposition \ref{delete-positive} again, (\ref{e3}) also holds for every $g_{i}$ with $i \in B$. However, note in passing that $g_{i}(b_{i}) = 0$ for all $i \in B$ since $b_{i}$ is not $\theta$-essential in $C_{i}$ (by Proposition \ref{Neumaier}).

Next, for each $i \in C$, set $g_{i}$ to be the zero vector on $V(C_{i})$. Note that (\ref{e3}) is satisfied vacuously for every $g_{i}$ with $i \in C$ since there are no $\theta$-essential vertices in the corresponding $C_{i}$.

Finally, we extend these $g_{i}$'s to an eigenvector of $T$ as follows: define $g \in \mathbb{R}^{|V(T)|}$ by
\[g(x) = \left\{\begin{array}{ll}

g_{i}(x)         &  ~~~~ \textrm{if } x \in V(C_{i})~~~\textrm{for some }i \in A \cup B \cup C,     \\
0      &   ~~~~\textrm{if } x = u.

\end{array} \right.   \]
Since (\ref{e3}) holds for every $g_{i}$,  we must have $g(x) \not = 0$ for every $\theta$-essential vertex $x$ of $T$. It is also readily verified that conditions in (\ref{eigenvalue}) are satisfied so that $g$ is indeed a $\theta$-eigenvector of $T$, as desired.

Moreover, by our construction, if $x$ is $\theta$-essential in $T \setminus u$, then $g(x) \not = 0$. By Proposition \ref{Neumaier}, $x$ must be $\theta$-essential in $T$, proving the second assertion of the theorem

\end{proof}

\noindent {\bf Proof of Theorem \ref{stability-tree}.}

Recall that $u$ is a given $\theta$-special vertex of $T$. By Proposition
\ref{delete-positive}, it remains to show that if $v$ is
$\theta$-essential in $T \setminus u$ then $v$ is $\theta$-essential
in $T$. But this is just the second assertion of the preceding theorem. This proves the Stability Lemma for trees.

\begin{example}\label{counter-example}
Let $T$ be the following tree \rm{:}
\end{example}
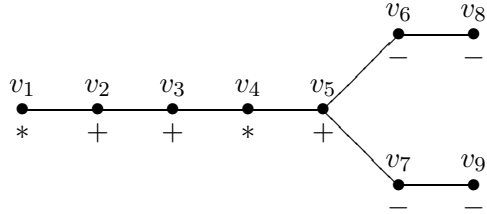
\[\xy
(0,0)*{}="A"; (10,0)*{}="B"; (20,0)*{}="C";
(30,0)*{}="D";(40,0)*{}="E"; (50,10)*{}="F"; (50,-10)*{}="G";
(60,10)*{}="H"; (60,-10)*{}="I";
"A";"E"**\dir{-};"E";"F"**\dir{-};"E";"G"**\dir{-};"F";"H"**\dir{-};"G";"I"**\dir{-};
(0,0)*{\bullet};
(10,0)*{\bullet};(20,0)*{\bullet};(30,0)*{\bullet};(40,0)*{\bullet};(50,10)*{\bullet};(50,-10)*{\bullet};
(60,10)*{\bullet};(60,-10)*{\bullet}; (0,3)*{v_1}="1";
(10,3)*{v_2}="v_2"; (20,3)*{v_3}="3";
(30,3)*{v_4}="4";(40,3)*{v_5}="5"; (50,13)*{v_6}="6";
(50,-7)*{v_7}="7"; (60,13)*{v_8}="8"; (60,-7)*{v_9}="9";
(0,-3)*{*}="1s"; (10,-3)*{+}="2s"; (20,-3)*{+}="3s";
(30,-3)*{*}="4s";(40,-3)*{+}="5s"; (50,7)*{-}="6s";
(50,-13)*{-}="7s"; (60,7)*{-}="8s"; (60,-13)*{-}="9s";
\endxy\]
The vertices are labeled $v_1, \ldots, v_9$ and the symbols $*$, $+$,
$-$ below each vertex indicates whether it is $\theta$-neutral or
$\theta$-positive or $\theta$-essential respectively where
$\theta=1$. Note that $\mu(T, x)=x^{9}-8x^{7}+20x^{5}-18x^{3}+5x$
and $\m(1, T)=1$. As the vertex $v_5$ is adjacent to a
$\theta$-essential vertex, $v_5$ is $\theta$-special in $T$. By
Theorem \ref{stability-tree}, upon deleting $v_{5}$ from $T$, all
other vertices are `stable' with respect to their $\theta$-partition
classes:

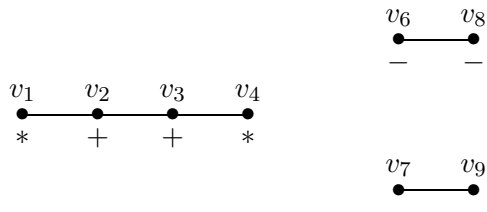
\[\xy
(0,0)*{}="A"; (10,0)*{}="B"; (20,0)*{}="C";
(30,0)*{}="D";(40,0)*{}="E"; (50,10)*{}="F"; (50,-10)*{}="G";
(60,10)*{}="H"; (60,-10)*{}="I";
"A";"D"**\dir{-};"F";"H"**\dir{-};"G";"I"**\dir{-}; (0,0)*{\bullet};
(10,0)*{\bullet};(20,0)*{\bullet};(30,0)*{\bullet};(50,10)*{\bullet};(50,-10)*{\bullet};
(60,10)*{\bullet};(60,-10)*{\bullet};(0,3)*{v_1}="1";
(10,3)*{v_2}="v_2"; (20,3)*{v_3}="3"; (30,3)*{v_4}="4";
(50,13)*{v_6}="6"; (50,-7)*{v_7}="7"; (60,13)*{v_8}="8";
(60,-7)*{v_9}="9"; (0,-3)*{*}="1s"; (10,-3)*{+}="2s";
(20,-3)*{+}="3s"; (30,-3)*{*}="4s"; (50,7)*{-}="6s";
(50,-13)*{-}="7s"; (60,7)*{-}="8s"; (60,-13)*{-}="9s";
\endxy\]
However, this is generally not true if we delete a non-special
vertex, for example, deleting $v_{3}$ from $T$ gives the following:

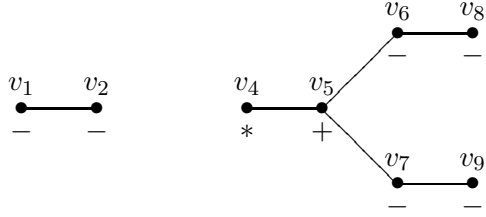
\[\xy
(0,0)*{}="A"; (10,0)*{}="B"; (20,0)*{}="C";
(30,0)*{}="D";(40,0)*{}="E"; (50,10)*{}="F"; (50,-10)*{}="G";
(60,10)*{}="H"; (60,-10)*{}="I";
"A";"B"**\dir{-};"D";"E"**\dir{-};"E";"F"**\dir{-};"E";"G"**\dir{-};"F";"H"**\dir{-};"G";"I"**\dir{-};
(0,0)*{\bullet};
(10,0)*{\bullet};(30,0)*{\bullet};(40,0)*{\bullet};(50,10)*{\bullet};(50,-10)*{\bullet};
(60,10)*{\bullet};(60,-10)*{\bullet};(0,3)*{v_1}="1";
(10,3)*{v_2}="2";
 (30,3)*{v_4}="4";(40,3)*{v_5}="5"; (50,13)*{v_6}="6";
(50,-7)*{v_7}="7"; (60,13)*{v_8}="8"; (60,-7)*{v_9}="9";
(0,-3)*{-}="1s"; (10,-3)*{-}="2s";
(30,-3)*{*}="4s";(40,-3)*{+}="5s"; (50,7)*{-}="6s";
(50,-13)*{-}="7s"; (60,7)*{-}="8s"; (60,-13)*{-}="9s";
\endxy\]

\section{Roots of Paths}

In this section, we prove some basic properties about roots of paths.

\begin{lm}\label{no-common}
Let $P_{n}$ denote the path on $n$ vertices, $n \ge 2$. Then $\mu(P_{n}, x)$ and $\mu(P_{n-1}, x)$ have no common root.
\end{lm}

\begin{proof}
Note that $\mu(P_{1}, x)=x$ and $\mu(P_{2}, x)=x^{2}-1$, and so they have no common root. Suppose $\mu(P_{n}, x)$ and $\mu(P_{n-1}, x)$ have a common root for some $n \ge 3$. Let $n$ be the least positive integer for which $\mu(P_{n},x)$ and $\mu(P_{n-1}, x)$ have a common root, say $\theta$. Then $\mu(P_{n-1},x)$ and $\mu(P_{n-2}, x)$ have no common root.

First we show that $\theta \not = 0$. Note that for any graph $G$, the multiplicity of $0$ as a root of its matching polynomial is the number of vertices missed by some maximum matching. Therefore, if $n$ is even then $P_{n}$ has a perfect matching, so $0$ cannot be a root of $\mu(P_{n}, x)$. It follows that if $n$ is odd then $0$ cannot be a root of $\mu(P_{n-1}, x)$. So $\theta \not = 0$.

Let $\{v_{1}, v_{2}\}$ be an edge in $P_{n}$ where $v_{1}$ is an endpoint of the path $P_{n}$. Note that $P_{n} \setminus v_{1} = P_{n-1}$ and $P_{n} \setminus v_{1}v_{2} = P_{n-2}$. By part (c) of Proposition \ref{identity}, $\mu(P_{n}, x) = x\mu(P_{n-1},x) - \mu(P_{n-2}, x)$, so $\theta$ is a root of $\mu(P_{n-1}, x)$ and $\mu(P_{n-2}, x)$, which is a contradiction. Hence $\mu(P_{n}, x)$ and $\mu(P_{n-1},x)$ have no common root.
\end{proof}

\begin{cor}\label{endpoints}
Let $\theta$ be a root of $\mu(P_{n},x)$. Then the endpoints of $P_{n}$ are $\theta$-essential.
\end{cor}

\begin{proof}
Suppose $v$ is an endpoint of $P_{n}$. If $v$ is $\theta$-neutral or $\theta$-positive in $P_{n}$ then $\theta$ is a root of $\mu(P_{n-1}, x)$, a contrary to Lemma \ref{no-common}.
\end{proof}

\begin{cor}\label{special}
Let $\theta$ be a root of $\mu(P_{n}, x)$. Then $P_{n}$ has no $\theta$-neutral vertices. Moreover, every $\theta$-positive vertex in $P_{n}$ is $\theta$-special.
\end{cor}

\begin{proof}
Let $v$ be a vertex of $P_{n}$ such that it is not $\theta$-essential. In view of Lemma \ref{neutral-essential}, it is enough to show that $v$ has a $\theta$-essential neighbor. By Corollary \ref{endpoints}, $v$ cannot be an endpoint of $P_{n}$. Then $P_{n} \setminus v$ consists of two disjoint paths, say $Q_{1}$ and $Q_{2}$. Let $u_{1}$ be the endpoint of $Q_{1}$ such that it is a neighbor of $v$ in $P_{n}$.

Consider the paths $Q_{1}$ and $Q_{1}v$ in $P_{n}$. Since $v$ is not $\theta$-essential, by Lemma \ref{essential-path}, $Q_{1}v$ is not $\theta$-essential in $P_{n}$. So the path $Q_{2} = P_{n} \setminus Q_{1}v$ has $\theta$ as a root of its matching polynomial.

If $Q_{1}$ is not $\theta$-essential in $P_{n}$ then the path $P_{n} \setminus Q_{1}$ would also have $\theta$ as a root of its matching polynomial. Since $P_{n} \setminus Q_{1}$ and $Q_{2}$ differ by exactly one vertex, this contradicts Lemma \ref{no-common}. Therefore, $Q_{1}$ is a $\theta$-essential path in $P_{n}$. By Lemma \ref{essential-path}, $u_{1}$ must be $\theta$-essential in $P_{n}$. Since $u_{1}$ is joined to $v$, we deduce from Lemma \ref{neutral-essential} that $v$ must be $\theta$-special.
\end{proof}

\section{Proof of Main Result}

We begin by proving the following special case.

\begin{prop}\label{m=2}
Let $T$ be a tree and $\m(\theta, T)=2$. Let $\mathcal{Q}=\{Q_{1}, Q_{2}\}$ be a set of vertex disjoint paths that cover $T$. Then $\mathcal{Q}$ is $(\theta, T)$-extremal.
\end{prop}

\begin{proof}
Since $T$ is a tree, there is an edge $\{u,v\} \in E(T)$ with $u \in V(Q_{1})$ and $v\in V(Q_{2})$. By Lemma \ref{interlacing-path}, $\m(\theta, Q_{1})=\m(\theta, T \setminus Q_{2}) \ge \m(\theta, T)-1 =1$. Similarly, $\m(\theta, Q_{2}) \ge 1$. Therefore $\theta$ is a root of $\mu(Q_{1},x)$ and $\mu(Q_{2}, x)$.

It remains to show that either $u$ is $\theta$-special in $Q_{1}$ or $v$ is $\theta$-special in $Q_{2}$. If all vertices in $T$ are $\theta$-essential, then $\m(\theta, T)=1$ by Theorem \ref{Neu}, which is impossible. So, there must be a $\theta$-special vertex in $T$, say $w$.

Suppose $w=u$. We shall prove that $w$ is also $\theta$-special in $Q_{1}$. Note that $\m(\theta, T \setminus w) = 3$. If $w$ is an endpoint of $Q_{1}$ then $T \setminus w$ is a disjoint union of two paths $Q_{1} \setminus w$ and $Q_{2}$. Since $Q_{1} \setminus w$ and $Q_{2}$ cover $T \setminus w$, we deduce from Theorem \ref{max-path} that $\m(\theta, T \setminus w) \le 2$, a contradiction. So $w$ is not an endpoint of $Q_{1}$. Removing $w$ from $Q_{1}$ would result in two disjoint paths, say $R_{1}$ and $R_{2}$. Note that $T \setminus w$ is the disjoint union of $R_{1}$, $R_{2}$ and $Q_{2}$. By part (a) of Proposition \ref{identity}, $\m(\theta, T \setminus w) = \m(\theta, R_{1})+\m(\theta, R_{2})+\m(\theta, Q_{2})$. By Theorem \ref{path} and the fact that $\m(\theta, T \setminus w)=3$, we conclude that $\m(\theta, R_{1})=\m(\theta, R_{2})=\m(\theta, Q_{2})=1$. Therefore, $\m(\theta, Q_{1} \setminus w) = \m(\theta, R_{1})+\m(\theta, R_{2})=2$. This means that $w$ must be $\theta$-positive in $Q_{1}$. By Corollary \ref{special}, $w$ is $\theta$-special in $Q_{1}$, as desired.

The case $w=v$ can be proved similarly.

Therefore, we may assume that $w \not = u, v$. We now proceed by induction on the number of vertices. Without loss of generality, we may assume that $w \in V(Q_{1})$. As before, it can be shown that $w$ is not an endpoint of $Q_{1}$.  So removing $w$ from $Q_{1}$ results in two disjoint paths, say $S_{1}$ and $S_{2}$. We may assume that $u \in V(S_{2})$. Then $T \setminus w$ is a disjoint union of $S_{1}$ and $T'$ where $T'$ is the tree induced by $S_{2}$ and $Q_{2}$. By Theorem \ref{path}, $\m(\theta, S_{1}) \le 1$. Since $S_{2}$ and $Q_{2}$ cover $T'$, by Theorem \ref{max-path}, $\m(\theta, T') \le 2$. As $w$ is $\theta$-special, $\m(\theta, T \setminus w)=3$. By part (a) of Proposition \ref{identity}, $\m(\theta, T \setminus w) = \m(\theta, S_{1}) + \m(\theta, T')$. We deduce that $\m(\theta, S_{1})=1$ and $\m(\theta, T')=2$. By induction, either $u$ is $\theta$-special in $S_{2}$ or $v$ is $\theta$-special in $Q_{2}$. In the latter, we are done. Therefore, we may assume that $u$ is $\theta$-special in $S_{2}$. So $u$ is not $\theta$-essential in $Q_{1} \setminus w$. Since $\m(\theta, Q_{1} \setminus w) = \m(\theta, S_{1}) + \m(\theta, S_{2})=2$, $w$ is $\theta$-positive in $Q_{1}$, so $w$ is $\theta$-special in $Q_{1}$ by Corollary \ref{special}. By the Stability Lemma for trees (Theorem \ref{stability-tree}), $u$ is not $\theta$-essential in $Q_{1}$. By Corollary \ref{special}, $u$ is $\theta$-special in $Q_{1}$.

Note that the base cases of our induction occur when $w=u$ or $w=v$.
\end{proof}

\begin{thm}\label{main1}
Let $T$ be a tree and $\m(\theta, T)=m$. Suppose $\mathcal{Q}=\{Q_{1}, \ldots, Q_{m}\}$ be a set of vertex disjoint paths that cover $T$. Then $\mathcal{Q}$ is $(\theta, T)$-extremal.
\end{thm}

\begin{proof}
We shall prove this by induction on $m \ge 1$. The theorem is trivial if $m=1$. If $m=2$, then the result follows from Proposition \ref{m=2}. So let $ \ge 3$. Since $T$ is a tree, there exist two paths, say $Q_{1}$ and $Q_{m}$, such that exactly one vertex in $Q_{1}$ is joined to other paths in $\mathcal{Q}$ and exactly one vertex in $Q_{m}$ is joined to other paths in $\mathcal{Q}$. To be precise, let $T'$ denote the tree induced by $Q_{2}, \ldots, Q_{m-1}$. Then there is only one edge joining $Q_{1}$ to $T'$ and only one edge joining $Q_{m}$ to $T'$.

By Theorem \ref{path}, $\m(\theta, T \setminus Q_{1}) \ge \m(\theta, T)-1=m-1$. Let $T''$ be the tree induced by $T'$ and $Q_{m}$, that is $T''=T \setminus Q_{1}$. Now $T''$ can be covered by $Q_{2}, \ldots, Q_{m}$. By Theorem \ref{max-path}, $\m(\theta, T'') \le m-1$. Therefore, $\m(\theta, T'') = m-1$. Moreover, $m-1$ is the maximum multiplicity of a root of $\mu(T'', x)$. By induction, $\{Q_{2}, \ldots, Q_{m}\}$ is $(\theta, T \setminus Q_{1})$-extremal.

 By a similar argument, $\{Q_{1}, \ldots, Q_{m-1}\}$ is $(\theta, T \setminus Q_{m})$-extremal. Hence $\mathcal{Q}$ is $(\theta, T)$-extremal.
\end{proof}

\begin{thm}\label{main2}
Let $F$ be a forest and $\mathcal{Q}=\{Q_{1}, \ldots, Q_{m}\}$ be a set of vertex disjoint paths that cover $F$. Suppose $\mathcal{Q}$ is $(\theta, F)$-extremal. Then $\m(\theta, F)=m$ and $\theta$ is a root of $\mu(F,x)$ with the maximum multiplicity.
\end{thm}

\begin{proof}
Since $F$ can be covered by $m$ vertex disjoint paths, by Theorem \ref{max-path}, we must have $\m(\alpha, F) \le m$ for any root $\alpha$ of $\mu(F,x)$. It remains to show that $\m(\theta, F) \ge m$.

An edge $\{u,v\}$ of $F$ is said to be $\mathcal{Q}$-{\em crossing} if $u$ and $v$ belong to different paths in $\mathcal{Q}$. If $F$ contains no $\mathcal{Q}$-crossing edges then $F$ consists of $m$ disjoint paths $Q_{1}, \ldots, Q_{m}$. Clearly, $\m(\theta, F) =\sum_{i=1}^{m} \m(\theta, Q_{i})= m$, as required. So we may assume that there exists an edge $\{u,v\} \in E(F)$ such that $u \in V(Q_{1})$ and $v \in V(Q_{2})$. Since $\mathcal{Q}$ is $(\theta, F)$-extremal, either $u$ is $\theta$-special in $Q_{1}$ or $v$ is $\theta$-special in $Q_{2}$.

We now proceed by induction on the number of vertices. Suppose $u$ is $\theta$-special in $Q_{1}$. Since $u$ is not an endpoint of $Q_{1}$, $Q_{1} \setminus u$ consists of two disjoint paths, say $R_{1}$ and $R_{2}$. Since $\m(\theta, Q_{1} \setminus u) = 2$ and $\m(\theta, R_{i}) \le 1$ for each $i=1,2$ (by Theorem \ref{path}), we deduce that $\m(\theta, R_{i})=1$ for each $i=1,2$. Note that $\{R_{1},R_{2}, Q_{3}, \ldots, Q_{m}\}$ is a set of disjoint paths that cover $F \setminus u$. Recall that there are no $\theta$-neutral vertices in $Q_{1}$. Moreover, by the Stability Lemma for trees (Theorem \ref{stability-tree}), every $\theta$-positive vertex in $Q_{1}$ remains $\theta$-positive in $Q_{1} \setminus u$ and every $\theta$-essential vertex in $Q_{1}$ remains $\theta$-essential in $Q_{1} \setminus u$. So every $\theta$-special vertex in $Q_{1}$ remains $\theta$-special in $Q_{1} \setminus u$. Consequently,  $\{R_{1}, R_{2}, Q_{3}, \ldots, Q_{m}\}$ is $(\theta, F \setminus u)$-extremal. By induction, $\m(\theta, F \setminus u)=m+1$ and $\theta$ is a root of $\mu(F \setminus u, x)$ with maximum multiplicity. It follows from Lemma \ref{interlacing} that $\m(\theta, F) \ge \m(\theta, F \setminus u)-1 = m$, as desired.

The case when $v$ is $\theta$-special in $Q_{2}$ can be settled by a similar argument. Note that the base cases of our induction occur when $F$ has no crossing edges.
\end{proof}

Our main result Theorem \ref{main} now follows immediately from Theorem \ref{main1} and Theorem \ref{main2}.

\end{document}